\input amstex
\input epsf
\magnification=\magstep1 
\baselineskip=13pt
\documentstyle{amsppt}
\def\rk{\operatorname{rank}}

\def\EE{\bold{E\thinspace}}
\def\var{\bold{var\thinspace}}
\def\spa{\operatorname{span}}
\def\Cov{\bold{Cov\thinspace}}
\def\vl{\operatorname{vol}}

\def\Es{\Cal{E}}

\vsize=8.7truein \CenteredTagsOnSplits \NoRunningHeads
\topmatter
 
\title A quick estimate for the volume of a polyhedron \endtitle 
\author Alexander Barvinok and Mark Rudelson  \endauthor
\address Department of Mathematics, University of Michigan, Ann Arbor,
MI 48109-1043, USA \endaddress
\email barvinok$\@$umich.edu, rudelson$\@$umich.edu \endemail
\date June 23,  2022 \enddate
\thanks  Research of AB and MR is partially supported by NSF Grants DMS 1855428 and DMS 2054408 respectively. 
\endthanks 
\abstract  Let $P$ be a bounded polyhedron defined as the intersection of the non-negative orthant ${\Bbb R}^n_+$ and an affine subspace of codimension $m$ in ${\Bbb R}^n$.
We show that a simple and computationally efficient formula approximates the volume of $P$ within a factor of $\gamma^m$, where $\gamma >0$ is an absolute constant. 
The formula provides the best known estimate for the volume of transportation polytopes from a wide family.
\endabstract
\endtopmatter
\document

\head 1. Introduction \endhead

The problem of efficient computation (approximation) of the volume of a polytope, and, more generally, of a given convex body has attracted a lot of attention, see, for example, 
\cite{GK18} for a survey. The most successful approach is via Markov Chain Monte Carlo randomized algorithms, see \cite{Ve05} for a survey. In particular, randomized algorithms allow one to approximate the volume of a polytope in ${\Bbb R}^n$ within relative error $\epsilon >0$ in time polynomial in $n$ and $\epsilon^{-1}$. The polytope can be defined as the convex hull of a finite set of points or as the intersection of halfspaces, or by a membership oracle, in which case the algorithms extend to the class of all ``well-conditioned" convex bodies.

Deterministic algorithms enjoyed less success. For a general convex body $B \subset {\Bbb R}^n$, the only available polynomial time algorithm approximates volume within a factor of $n^{O(n)}$, using an approximation of $B$ by an ellipsoid, see \cite{G+88}. For $B$ defined by a membership oracle, this approximation factor is basically the best possible (up to some logarithmic terms) that can be achieved in deterministic polynomial time. More precisely, it is shown in \cite{BF87} that with $N$ queries to the membership oracle, one cannot deterministically estimate the volume better than within a factor of $\left(\gamma n/ \ln N \right)^{n/2}$, where $\gamma >0$ is an absolute constant. A deterministic algorithm approximating the volume within a factor of $2^{O(n)}$ in $2^{O(n)}$ time, which matches the above mentioned lower bound, is constructed in \cite{DV13}. For a convex body 
$B \subset {\Bbb R}^n$ defined by a membership oracle, for any $0< \epsilon <1$, the algorithm from \cite{Da15}, see also \cite{DV13}, approximates $\vl B$ within a factor of $(1+ \epsilon)^n$ in 
$O\left((1/\epsilon)^{O(n)}\right)$ time, almost matching a lower bound that follows from \cite{BF88}, see also \cite{B+89}, \cite{CP88}, and \cite{Gl89}.

If $P$ is a polytope defined as the convex hull of a set of points or as the intersections of halfspaces, deterministic algorithms in principle may turn out to be as powerful as randomized ones, but so far the approximation ratio achieved in deterministic polynomial time is the same as for general convex bodies. We remark that if $P \subset {\Bbb R}^n$ is a polytope defined as the convex hull of $n +O(1)$ points or as the intersection of $n +O(1)$ halfspaces, then $\vl P$ can be computed exactly in polynomial time, in the former case by a triangulation into $n^{O(1)}$ simplices and in the latter case by a dual procedure of expressing $P$ as a signed linear combination of $n^{O(1)}$ simplices, see \cite{GK18} and \cite{La91}. If $P$ is defined by a system of linear inequalities with rational coefficients, then even to write $\vl P$ as a rational number, one may need the number of bits that is exponential in the size of the input \cite{La91}. For polytopes $P$ defined by a system of linear inequalities with a totally unimodular matrix of integer coefficients, as well as for polytopes defined as convex hulls of sets of rational points,  the problem of computing $\vl P$ exactly is $\#$P-hard \cite{DF88}.
Of course, for special classes of polytopes, such as parallelepipeds, there can be computationally efficient explicit formulas.

In this paper, we consider the class of polyhedra $P$ defined as the intersection of the non-negative orthant ${\Bbb R}^n_+$ and an affine subspace in ${\Bbb R}^n$. In coordinates, $P$ is defined by a system of linear equations $Ax=b$, where $A$ is an $m \times n$ matrix, $x$ is an $n$-vector of variables and $b$ is an $m$-vector, and inequalities $x \geq 0$, meaning that the coordinates of $x$ are non-negative. We assume that $m <n$, that $\rk A=m$ and that $P$ has a non-empty relative interior, that is, contains a point $x$ with all coordinates positive. Hence 
$\dim P=n-m$ and we measure the $(n-m)$-dimensional volume of $P$ in its affine span with respect to the Euclidean structure inherited from ${\Bbb R}^n$. We also assume that $P$ is bounded, that is, a polytope. Generally, any $(n-m)$-dimensional polyhedron with $n$ facets can be represented as the intersection of ${\Bbb R}^n_+$ and an affine subspace of codimension $m$. Furthermore, many interesting polyhedra, such as transportation polytopes, see, for example \cite{DK14}, are naturally defined in this way. 

We present a deterministic polynomial time algorithm which approximates the volume of such a polytope $P$ within a factor of $\gamma^m$, where $\gamma >0$ is an absolute constant (for $m$ large enough, one can choose $\gamma=4.89$). In fact, our algorithm is basically a {\it formula}. The only ``non-formulaic" part of our algorithm consists of solving some standard convex optimization problem on $P$, namely finding its ``analytic center", see \cite{Re88}. 
After that, we only need to compute two $m \times m$ determinants, which, as is well-known, can be accomplished  in 
$O(m^3)$ time.
While the approximation factor $\gamma^m$ looks big compared to $1+\epsilon$ achieved by randomized algorithms, it appears to be the best achieved to date by a deterministic polynomial time algorithm for many interesting classes of polytopes, such as transportation polytopes. Since the algorithm is basically a formula, it allows one to analyze how the volume changes as $P$ evolves inside its class, which turns out to be important for studying some statistical phenomena related to contingency tables, cf. \cite{D+20}.
The approximation factor looks more impressive when $n \gg m$, which is indeed the case for many interesting classes of polytopes. Note that if we dilate a $d$-dimensional polytope by a factor of $(1+\epsilon)$, its volume gets multiplied by $(1+\epsilon)^d$. 
This has implications for evaluating the volume ratio defined by 
$$v(P)=\left({\vl P \over \vl {\Cal B}}\right)^{1/d},$$
where ${\Cal B}$ is the maximum volume ellipsoid inscribed in $P$. This quantity plays a fundamental role in geometric functional analysis \cite{Pi89}, \cite{A+15}. As the volume of ${\Cal B}$ can be efficiently calculated, see \cite{GK18} and references therein, we can approximate the volume ratio within a factor of $1+o(1)$ in deterministic polynomial time when $n \gg m$.

\head 2. The main result and some applications \endhead

\subhead (2.1) The setup \endsubhead
Let $A=\left(\alpha_{ij}\right)$ be an $m \times n$ matrix, let $b=\left(\beta_1, \ldots, \beta_m\right)$ be an $m$-vector, and suppose that the polyhedron $P \subset {\Bbb R}^n$ 
is defined by the system of equations 
$$\sum_{j=1}^n \alpha_{ij} \xi_j =\beta_i \quad \text{for} \quad i=1, \ldots, m \tag2.1.1$$
and inequalities 
$$\xi_j \geq 0 \quad \text{for} \quad j=1, \ldots, n. \tag2.1.2$$
We assume that $m < n$ and that $\rk A =m$, in which case the system (2.1.1) defines an $(n-m)$-dimensional affine subspace.

Suppose that $P$ has a non-empty relative interior, that is, contains a point $x=\left(\xi_1, \ldots, \xi_n\right)$ where $\xi_j > 0$ for $j=1, \ldots, n$, and is also bounded, that is, 
a polytope. Thus $P$ is an $(n-m)$-dimensional bounded polyhedron and our goal is to estimate its volume $\vl P$ relative to the Euclidean structure inherited from ${\Bbb R}^n$ 
by the affine subspace defined by (2.1.1).

We define a function $f: {\Bbb R}^n_+ \longrightarrow {\Bbb R}$ by 
$$\split &f(x)=n + \sum_{j=1}^n \ln \xi_j \quad \text{where} \\& x=\left(\xi_1, \ldots, \xi_n\right) \quad \text{and} \quad \xi_j > 0 \quad \text{for} \quad j=1, \ldots, n 
\endsplit \tag2.1.3$$
and consider the following optimization problem
$$\text{Find} \quad z \in P \quad \text{such that} \quad f(z)=\max_{x \in P} f(x). \tag2.1.4$$
The function $f$ is strictly concave and hence the maximum point $z$ can be found efficiently (in polynomial time), see \cite{NN94}. 
Also, the point $z=\left(\zeta_1, \ldots, \zeta_n\right)$ is unique and satisfies $\zeta_j >0$ for $j=1, \ldots, n$, see \cite{BH10}. In fact, the point $z$ was at the very source of interior-point methods in optimization \cite{Re88} under the name of the {\it analytic center} of $P$. We note that if the symmetry group of $P$ is sufficiently rich, we can determine $z$ without solving the optimization problem (2.1.3)--(2.1.4), as $z$ must be invariant under all permutations of the coordinates $\xi_1, \ldots, \xi_n$ that map $P$ onto itself.

Now we can state our result.
\proclaim{(2.2) Theorem} Let $A$ be $m \times n$ matrix of $\rk A =m < n$, let $b$ be an $m$-vector and suppose that the polyhedron $P$ defined by the system (2.1.1)--(2.1.2) is bounded and has a non-empty relative interior. Let $z=\left(\zeta_1, \ldots, \zeta_n\right)$ be the analytic center of $P$ defined as the solution to the optimization problem (2.1.3)--(2.1.4). 
Let $B$ be the $m \times n$ matrix obtained by multiplying the $j$-th column of $A$ by $\zeta_j$ for $j=1, \ldots, n$ and let 
$$\Es(A, b) = e^{f(z)} {\sqrt{\det A A^T} \over \sqrt{\det B B^T}}=e^n \zeta_1 \cdots \zeta_n {\sqrt{\det A A^T} \over \sqrt{\det B B^T}}.$$
\roster
\item Let $\alpha_0$ be the necessarily unique number in the interval $(0, 1)$ satisfying
$${1 \over 2\pi} \int_{-\infty}^{+\infty} \left(1+ s^2\right)^{-{1 \over 2\alpha_0}}\ ds =1, \quad \alpha_0 \approx 0.7148659168.$$
Then
$$\vl P \ \leq \ \left({1 \over \sqrt{\alpha_0}}\right)^m \Es(A,b) \ \leq \ (1.19)^m \Es(A,b);$$
\item We have 
$$\vl P \ \geq \ {2 \Gamma\left({m +2 \over 2}\right) \over \pi^{m/2} e^{(m+2)/2} (m+2)^{m/2}} \Es(A, b);$$
\item For any $0 < \epsilon < 1/2$, there is $\gamma(\epsilon) >0$ such that 
$$\split \vl P \ \geq \ &\exp\left\{ -\gamma(\epsilon) \sqrt{m} \ln^4 (m+1) \right\} \left({1 \over \sqrt{2 \pi e}}\right)^m \Es(A, b) \\
\ \geq \ & \exp\left\{ -\gamma(\epsilon) \sqrt{m} \ln^4 (m+1) \right\}  (0.24)^m \Es(A, b). \endsplit $$
\endroster
\endproclaim

Some remarks are in order. Using the standard bound 
$$\Gamma(t) \ \geq \ \sqrt{2 \pi} t^{t-{1 \over 2}} e^{-t} \quad \text{for} \quad t \geq 1,\tag2.2.1$$
we conclude that the right hand side of the formula in Part (2) decreases with $m$ roughly as 
$$\left({1 \over e \sqrt{2 \pi}}\right)^m \Es(A, b) \approx (0.14)^m \Es(A, b).$$ The lower bound in Part (3) is asymptotically stronger, although it contains a constant which may render it weaker than the bound of Part (2) for smaller $m$. In addition, the proof of Part (2) is rather elementary, whereas the proof of Part (3) relies on the recent breakthrough by Chen \cite{Ch21} and  Klartag and Lehec \cite{KL22} in thin shell estimates,
although the same asymptotic behavior in the $\left({1\over \sqrt{2 \pi e}}\right)^m$ term could be achieved by using earlier and weaker, but still highly non-trivial bounds from \cite{Kl07}.

We also note that the estimate $\Es(A, b)$ scales properly when the polyhedron is dilated: if $b\longmapsto \tau b$ for $\tau >0$, then $z \longmapsto \tau z$ and 
$\Es(A, \tau b) =\tau^{n-m} \Es(A, b)$.

\subhead  (2.3) Gaussian approximation \endsubhead 
It was proved in \cite{BH10} that if some analytic conditions on $A$ and $z$ are satisfied, we have asymptotically
$$\vl P \approx {e^{f(z)} \over (2\pi)^{m/2}} {\sqrt{\det A A^T} \over \sqrt{\det B B^T}}={e^n \zeta_1 \cdots \zeta_n \over (2\pi)^{m/2}}  {\sqrt{\det A A^T} \over \sqrt{\det B B^T}}  \tag2.3.1$$
as $m$ and $n$ grow. The right hand side of (2.3.1) is called in \cite{BH10} the {\it maximum entropy Gaussian approximation}. Under typical circumstances those analytic conditions require, in particular, that $m=O(\sqrt{n})$ and that the coordinates $\zeta_1, \ldots, \zeta_n$ of the analytic center $z$ of $P$ are roughly of the same order.
We explain the name and the intuition behind this formula in Section 3.1. The estimate of Theorem 2.2 is much cruder, but its validity doesn't depend on the particulars of $A$ and $b$ or the relations between $m$ and $n$. We note that to understand some statistical phenomena related to contingency tables \cite{D+20}, it is important to understand the behavior of the volume of $P$ when the coordinates $\zeta_1, \ldots, \zeta_n$ have decidedly different orders of magnitude.

\subhead (2.4) Example: simplex \endsubhead Suppose that $m=1$, so $P$ is defined by a single linear equation 
$$\alpha_1 \xi_1 + \ldots + \alpha_n \xi_n =\beta$$ together with the inequalities
$$\xi_j \geq 0 \quad \text{for} \quad j=1, \ldots, n.$$ 
Since we assume that $P$ is bounded and with a non-empty relative interior, we must have $\beta >0$ and $\alpha_j >0$ for $j=1, \ldots, n$. Since $\Es(A, b)$ scales correctly when $b$ is scaled, 
$b \longrightarrow \tau b$ for $\tau > 0$,  we further assume
that 
$\beta=n$.
Then for the analytic center $z=\left(\zeta_1, \ldots, \zeta_n\right)$, we have 
 $$\zeta_j = {1 \over \alpha_j} \quad \text{for} \quad j=1, \ldots, n \tag2.4.1$$
 and hence
 $$\Es(A, b)= {e^n \sqrt{\alpha_1^2 + \ldots + \alpha_n^2} \over \alpha_1 \cdots \alpha_n \sqrt{n}}.$$
 On the other hand, 
 $$\vl P= {n^n \sqrt{\alpha_1^2 + \ldots + \alpha_n^2} \over n! \alpha_1 \cdots \alpha_n}.$$
 By Stirling's formula
 $$n! = n^n e^{-n} \sqrt{2\pi n} \left(1+ o(1)\right) \quad \text{as} \quad n \longrightarrow \infty,$$
 and hence the Gaussian approximation is asymptotically exact as $n \longrightarrow \infty$.

\subhead (2.5) Example: $2$-way transportation polytopes \endsubhead Let us fix positive integers $k$ and $l$, a $k$-vector of positive real numbers $r=\left(\rho_1, \ldots, \rho_k\right)$ and 
an $l$-vector of positive real numbers $c=\left(\gamma_1, \ldots, \gamma_l\right)$ such that 
$$\sum_{i=1}^k \rho_i = \sum_{j=1}^l \gamma_j. \tag2.5.1$$
We consider the set $T(r, c)$ of $k \times l$ non-negative real matrices with row sums $r$ and column sums $c$. In other words, $T(r, c)$ is defined in the space 
${\Bbb R}^{k \times l} \cong {\Bbb R}^{kl}$ of $k \times l$ matrices $x=\left(\xi_{ij}\right)$ by the equations 
$$ \sum_{j=1}^l \xi_{ij} = \rho_i \quad \text{for} \quad i=1, \ldots, k \quad \text{and} \quad \sum_{i=1}^k \xi_{ij} = \gamma_j \quad \text{for} \quad j=1, \ldots, l \tag2.5.2$$
and inequalities 
$$\xi_{ij} \ \geq \ 0 \quad \text{for all} \quad i, j.$$
It is not hard to see that $T(r, c)$ is a polyhedron of dimension $(k-1)(l-1)$. Because of the balance condition (2.5.1), which is necessary and sufficient for $T(r, c)$ to be non-empty, the 
equations (2.5.2) are not linearly independent, and to bring them into the form required by Theorem 2.2 it suffices to drop precisely one equation from the list. The polyhedron $T(r, c)$ is called a {\it 2-way transportation polytope} with {\it margins} $r$ and $c$, see \cite{DK14}, and its volume was investigated, in particular, in \cite{CM09}, \cite{Ba09}, \cite{BH12} and
\cite{B+20}.

By symmetry, it follows that for the analytic center $z=\left(\zeta_{ij}\right)$ we have 
$$\zeta_{ij}={\rho_i \over l} \quad \text{provided} \quad \gamma_1=\ldots = \gamma_l$$
and 
$$\zeta_{ij}={\gamma_j \over k} \quad \text{provided} \quad \rho_1=\ldots = \rho_k.$$

Some margins are of a particular interest. If $k=l$ and $\rho_i=\gamma_j=1$ for all $i$ and $j$, which we write as $r=c={\bold 1}$, 
we get the polytope of $k \times k$ doubly stochastic matrices, also known as the {\it Birkhoff} or {\it Birkhoff - von Neumann} polytope, cf. \cite{DK14}. In this case, Canfield and McKay \cite{CM09} obtained an asymptotic formula for the volume as $k \longrightarrow \infty$:
$$\vl T({\bold 1}, {\bold 1}) = {1 \over (2\pi)^{k -{1 \over 2}} k^{(k-1)^2}} \exp\left\{ {1 \over 3} + k^2\right\} \left(1 +o(1)\right).$$
By symmetry, the analytic center of $T({\bold 1}, {\bold 1})$ is the matrix $Z=\left(\zeta_{ij}\right)$ with 
$$\zeta_{ij}={1 \over k} \quad \text{for all} \quad i, j.$$
Note that the formula differs from the Gaussian approximation (2.3.1) by a factor of $e^{1/3}$. In \cite{BH12}, this factor was interpreted as the {\it Edgeworth correction} in the Central Limit Theorem, and so corrected Gaussian approximation asymptotic formula was extended for all ``tame" margins, where $k$ and $l$ grow proportionately, and all coordinates 
$\zeta_{ij}$ of the analytic center are within a constant factor of each other. 

One can observe a curious phase transition destroying the tameness of margins somewhat unexpectedly. Suppose that $k=l$ and that 
$$\rho_1=\ldots=\rho_{k-1}=\gamma_1= \ldots = \gamma_{k-1}=1.$$
It is not hard to show that if we choose 
$$\rho_k = \gamma_k = 2-\epsilon \quad \text{for some small} \quad \epsilon > 0$$
then the entries of the analytic center satisfy
$$\max_{ij} \zeta_{ij} =O\left(k^{-1}\right) \quad \text{as} \quad k \longrightarrow \infty.$$
However, if we choose 
$$\rho_k=\gamma_k=2+\epsilon \quad \text{for some small} \quad \epsilon >0,$$ then the $\zeta_{kk}$ entry becomes large:
$$\zeta_{kk} \ > \ \delta \quad \text{for some} \quad \delta =\delta(\epsilon)>0.$$
This and similar phase transitions are investigated in \cite{D+20}. Their existence may serve as an indication that $\vl T(r, c)$ cannot be estimated too closely by a smooth analytic expression as the margins $r$ and $c$ vary even mildly. The formula of Theorem 2.2 approximates $\vl T(r, c)$ within a factor of 
$\exp\left\{ O(k+l)\right\}$ and it appears to be the only known formula where the bound on the approximation factor does not depend on the numerics of the margins $r$ and $c$.
It also provides the best known approximation for the volume of a generic 2-way transportation polytope.

\subhead (2.6) Example: 3-way planar transportation polytopes \endsubhead
For an integer $r >0$ we consider the polytope $P_r$ of all $r \times r \times r$ arrays (tensors) $X=\left(\xi_{ijk}\right)$ satisfying the equations
$$\aligned &\sum_{i=1}^r \xi_{ijk} =1 \quad \text{for} \quad j,k=1, \ldots, r, \quad \sum_{j=1}^r \xi_{ijk}=1 \quad \text{for} \quad i,k=1, \ldots, r \\
&\text{and} \quad \sum_{k=1}^r \xi_{ijk} =1 \quad \text{for} \quad i,j=1, \ldots, r\endaligned \tag2.6.1$$
and inequalities
$$\xi_{ijk} \ \geq \ 0 \quad \text{for all}  \quad i, j, k.$$
The polytope $P_r$ is known as a {\it 3-way planar transportation polytope}, see \cite{DK14}. One can also consider  {\it 3-way axial transportation polytopes} obtained by fixing sums over 2-dimensional coordinate sections of the array; somehow, those turn out to have a simpler structure than $P_r$. The linear equations (2.6.1) are not independent, and it is not hard to check that $\dim P_r=(r-1)^3$.

As is well known, the vertices of the Birkhoff polytope of Section 2.5 are the permutation matrices. The {\it integer} vertices of $P_r$ correspond to Latin squares, but there are plenty of non-integer vertices \cite{LL14} and the arithmetic of their coordinates can vary wildly \cite{Gr92}. By symmetry, the analytic center $Z=\left(\zeta_{ijk}\right)$ of $P_r$ satisfies 
$$\zeta_{ijk}={1 \over r} \quad \text{for all} \quad i,j,k.$$
Theorem 2.2 implies that up to a factor of $\gamma^{r^2}$ for some absolute constant $\gamma >0$, the volume of $P_r$ is approximated by 
$e^{r^3} r^{-(r-1)^3}$. Hence we obtain an asymptotically exact estimate
$$\ln \vl P_r = r^3 -(r-1)^3 \ln r + O(r^2) \quad \text{as} \quad r \longrightarrow \infty. \tag2.6.2$$
It appears that (2.6.2) is the best estimate of the volume of $P_r$ to date.

In the rest of the paper, we prove Theorem 2.2. In Section 3, we collect some preliminaries. In Section 4, we prove the upper bound of Part (1). In Section 5, we prove the lower bound of Part (2), and in Section 6, we prove the lower bound of Part (3).

\head 3. Preliminaries \endhead

\subhead (3.1) The maximum entropy density \endsubhead Recall that a real-valued random variable $X$ has the {\it standard exponential distribution} if the density 
$p_X(t)$ of $X$ satisfies 
$$p_X(t)=\cases e^{-t} &\text{if\ } t \geq 0 \\ 0 &\text{if\ } t < 0. \endcases$$
For the expectation and variance, we have
$$\EE X=1 \quad \text{and} \quad \var X=1.$$
Let matrix $A$, vector $b$, function $f$, polyhedron $P$ and point $z=\left(\zeta_1, \ldots, \zeta_n\right)$ be as in Theorem 2.2. Let
$a_1, \ldots, a_n$ be the columns of matrix $A$, considered as $m$-vectors. Suppose further that $X_1, \ldots, X_n$ are independent standard exponential random variables 
and let us define a random vector $Y$ with values in ${\Bbb R}^m$ by 
$$Y=\sum_{j=1}^n \zeta_j X_j a_j = \sum_{j=1}^n X_j b_j, \tag3.1.1$$
where $b_1, \ldots, b_n$ are the columns of matrix $B$, as defined in Theorem 2.2.
It is proved in \cite{BH10} that the density $p_Y$ at $b=\left(\beta_1, \ldots, \beta_m\right)$ can be expressed as 
$$p_Y(b)={\vl P \over e^{f(z)} \sqrt{\det A A^T}}\tag3.1.2$$
and that for the expectation and the covariance matrix of $Y$, we have 
$$\EE Y=b \quad \text{and} \quad \Cov Y=BB^T. \tag3.1.3$$
We also need the characteristic function of $Y$. For $t\in {\Bbb R}^m$, $t=\left(\tau_1, \ldots, \tau_m\right)$, we have
$$\split \phi_Y(t) = &\EE \exp\Bigl\{ \sqrt{-1} \langle Y, t \rangle \Bigr\} =\prod_{j=1}^n \EE \exp\Bigl\{ \sqrt{-1}  X_j \langle b_j, t\rangle\Bigr\}\\
= &\prod_{j=1}^n {1 \over 1- \sqrt{-1} \langle b_j, t \rangle},\endsplit$$
where $\langle \cdot, \cdot \rangle$ is the standard inner product in ${\Bbb R}^m$. Consequently, the density of $Y$ can be recovered 
as 
$$p_Y(y)={1 \over (2 \pi)^m} \int_{{\Bbb R}^m} \exp\Bigl\{-\sqrt{-1} \langle y, t \rangle \Bigr\} \prod_{j=1}^n {1 \over 1- \sqrt{-1} \langle b_j, t \rangle} \ dt, \tag3.1.4$$
see also \cite{BH10}. We will be interested only in the situations when the integral (3.1.4) converges absolutely.

Equations (3.1.1) -- (3.1.4) are the only ones we need from this section for the proof of Theorem 2.2. The rest contains some explanatory remarks.

Equations (3.1.3) and (3.1.4) are straightforward to check, while equation (3.1.2) follows from the fact that the density of the random vector $Z=\left(\zeta_1 X_1, \ldots, \zeta_n X_n \right)$ is constant on $P$ and equal to $e^{-f(z)}$, see Theorem 7 in \cite{BH10}. The formal proof easily follows from the Lagrange optimality condition for $z$. A more intuitive explanation is that $Z$ has the largest entropy among all random vectors supported on ${\Bbb R}^n_+$ and with expectation in the affine subspace defined by the system $Ax=b$, just as the standard exponential distribution has the largest entropy among all distributions supported on ${\Bbb R}_+$ and with expectation 1.

Since $Y$ is the sum of independent random variables and (3.1.3) holds, in view of the (local) Central Limit Theorem it is not inconceivable that in the vicinity of $b$, the distribution of $Y$ can be close to a Gaussian distribution. By analyzing the integral (3.1.4), it is shown in \cite{BH10} that it is indeed the case under some conditions on $A$ and $z$, and hence we obtain the Gaussian 
approximation formula (2.3.1). It is further shown in \cite{BH12} that for transportation polytopes with ``tame" margins, the local Central Limit Theorem holds, albeit with the Edgeworth correction that takes into account the 3rd and 4th moments of $Y$.

\subhead (3.2) Example: simplex \endsubhead Let $P$ be the simplex of Example 2.4, defined by the equation 
$$\alpha_1 \xi_1  + \ldots + \alpha_n \xi_n=n$$ 
and inequalities 
$$\xi_j \geq 0 \quad \text{for} \quad j=1, \ldots, n.$$ Then by (2.4.1) for the random variable $Y$ defined by (3.1.1), we have 
$$Y= \sum_{j=1}^n \zeta_j X_j \alpha_j = \sum_{j=1}^n X_j$$
and hence $Y$ is a random variable with gamma-density
$$p_Y(t)= \cases {t^{n-1}\over (n-1)!}  e^{-t} &\text{for\ } t \geq 0 \\ 0 &\text{for\ } t <0. \endcases$$
We have 
$$b=\EE Y=n \quad \text{and} \quad  f(z)= n - \sum_{j=1}^n \ln \alpha_j, $$
while the formula (3.1.2) reads 
$$\split p_Y(b)=&{n^{n-1} \over (n-1)! e^{n}}={n^n \sqrt{\alpha_1 + \ldots + \alpha_n^2} \over n!  \alpha_1 \cdots \alpha_n} \cdot { \alpha_1 \cdots \alpha_n \over  e^n\sqrt{\alpha_1^2 + \ldots +
\alpha_n^2}} \\=&{\vl P \over e^{f(z)} \sqrt{\det A A^T}}. \endsplit$$

\subhead (3.3) Isotropic and log-concave densities \endsubhead Recall that a non-negative measurable function $\rho: {\Bbb R}^m \longrightarrow {\Bbb R}_+$ is called 
{\it density} if 
$$\int_{{\Bbb R}^m} \rho(x) \ dx =1.$$
A density $\rho$ is called {\it centered} if 
$$\int_{{\Bbb R}^m} x \rho(x) \ dx =0,$$
where we assume that the integral converges absolutely.  A density $\rho$ is called {\it isotropic} if it is centered and
$$\int_{{\Bbb R}^m} \xi_i \xi_j \rho(x)\ dx =\cases 1 &\text{if\ } i=j \\ 0 &\text{if\ } i \ne j, \endcases \quad \text{where} \quad x=\left(\xi_1, \ldots, \xi_m\right).$$
A density $\rho: {\Bbb R}^m \longrightarrow {\Bbb R}_+$ is called {\it logarithmically concave} or {\it log-concave}, if it can be written as 
$\rho(x)=e^{\psi(x)}$, where $\psi: {\Bbb R}^m \longrightarrow {\Bbb R} \cup \{-\infty\}$ is a concave function. 

We will use the following basic fact.

Let
$\rho: {\Bbb R}^n \longrightarrow {\Bbb R}_+$ be a log-concave density, let $L \subset {\Bbb R}^n$ be a subspace and let ${\Bbb R}^n \longrightarrow L$ be the orthogonal projection. Then the push-forward density $\varphi$ on $L$ defined by 
$$\varphi(y) = \int_{y + L^{\bot}} \rho(x) \ dx,$$
where $L^{\bot}$ is the orthogonal complement of $L$, is 
also log-concave. This is a standard corollary of the Pr\'ekopa - Leindler inequality, see, for example, Chapter I of \cite{A+15}.

\subhead (3.4) Transforming the density of $Y$ \endsubhead
From (3.1.2), we express the volume of $P$ as 
$$\vl P = e^{f(z)} \sqrt{\det A A^T} p_Y(b), \tag3.4.1$$
where $p_Y$ is the density of the random variable $Y$ defined by (3.1.1).
The analytic center $z$ depends only on the affine subspace defined by the system $Ax=b$, but not on a particular choice of a matrix $A$ and vector $b$. 
If $W$ is an invertible $m \times m$ matrix, then the affine subspaces defined by systems $Ax=b$ and $A'x=b'$ where $A'=WA$ and $b'=Wb$ coincide. If we replace $A$ by $A'=WA$ 
and $b$ by $b'=Wb$ then the matrix $B$ gets replaced by $B'=WB$, and we have 
$${\det AA^T \over \det BB^T} = {\det A' (A')^T \over \det B' (B')^T},$$
and hence the estimate $\Es(A, b)$ of Theorem 2.2 does not change.
Furthermore, we have $B' (B')^T=W (B B^T) W^T$. Choosing an appropriate $W$ if needed, without loss of generality, we assume that $B$ satisfies 
$$B B^T = I_m, \tag3.4.2$$
where $I_m$ is the $m \times m$ identity matrix.
If (3.4.2) holds, then in view of the formula (3.1.3), the density $p_{Y-b}$ of $Y-b$ is isotropic. The crucial fact for us is that $p_Y$ is also log-concave.
Because of (3.4.2), we can identify ${\Bbb R}^m$ isometrically with an $m$-dimensional subspace $L$ in ${\Bbb R}^n$, so that $B$ is the matrix of the orthogonal projection 
${\Bbb R}^n \longrightarrow L$ in some pair of orthonormal bases of ${\Bbb R}^n$ and $L$. Let
$$\rho\left(\xi_1, \ldots, \xi_n\right)= \cases \exp\left\{ -\sum_{j=1}^n \xi_j \right\}   &\text{if\ } \xi_1, \ldots, \xi_n \geq 0 \\ 0 & \text{otherwise.} \endcases $$
be the standard exponential density on ${\Bbb R}^n$. Obviously, $\rho$ is log-concave.
Then $p_Y$ is the push-forward of $\rho$ and hence is also log-concave.

\head 4. Proof of the upper bound \endhead

In this section, we prove Part (1) of Theorem 2.2. Our approach is inspired by Ball's work on the volume of a section of the cube \cite{Ba89}.
\subhead (4.1) More preliminaries \endsubhead 
As is discussed in Section 3.4, we assume that matrix $B$ satisfies (3.4.2). We define the random variable $Y$ by (3.1.1).
In view of (3.4.1), 
our goal is to bound the $\ell^{\infty}$-norm $\left\|p_Y\right\|_{\infty}$ of the density $p_Y$ of $Y$. From (3.1.4), we have 
$$\left\|p_Y\right\|_{\infty} \ \leq \ {1 \over (2\pi)^m} \int_{{\Bbb R}^m} \prod_{j=1}^n \left( 1+ \langle b_j, t \rangle^2\right)^{-{1 \over 2}} \ dt \tag 4.1.1$$
For a vector $a \in {\Bbb R}^m$, $a=\left(\alpha_1, \ldots, \alpha_m\right)$, by 
$a \otimes a$ we denote the $m \times m$ matrix with the $(i, j)$-th entry equal $\alpha_i \alpha_j$. 
To bound the integral in the right hand side of (4.1.1), we will use the Brascamp - Lieb inequality in the form adapted by Ball, see Theorem 2 in \cite{Ba01}. 
\proclaim{(4.2) Lemma} Let $u_1, \ldots, u_n$ be unit vectors from ${\Bbb R}^m$ and let $\lambda_1, \ldots, \lambda_n$ be positive numbers such that 
$$\sum_{j=1}^n \lambda_j \left( u_j \otimes u_j \right)=I_m,$$
where $I_m$ is the $m \times m$ identity matrix. Then, for measurable functions $f_1, \ldots, f_n: {\Bbb R} \longrightarrow {\Bbb R}_+$, we have 
$$\int_{{\Bbb R}^m} \prod_{j=1}^n f_j^{\lambda_j} \left(\langle u_j, x \rangle\right) \ dx \ \leq \ \prod_{j=1}^n \left(\int_{-\infty}^{+\infty} f_j(\xi) \ d \xi \right)^{\lambda_j}.$$
\endproclaim
{\hfill \hfill \hfill} \qed

Next, we investigate 1-dimensional integrals.

\proclaim{(4.3) Lemma} For $\alpha \in (0, 1)$, let 
$$F(\alpha)={1 \over 2 \pi} \int_{-\infty}^{+\infty} \left( 1 + \alpha\tau^2\right)^{-{1 \over 2\alpha}} \ d \tau.$$
Then 
\roster
\item The function $F(\alpha)$ is increasing on the interval $(0, 1)$;
\item There is a unique $\alpha_0 \in (0, 1)$ such that 
$$F(\alpha_0)= {1\over \sqrt{\alpha_0}}.$$
\endroster
Numerically,
$$\alpha_0 \approx 0.7148659168.$$
\endproclaim
\demo{Proof}
Let 
$$h(\tau, \alpha)= \left( 1 + \alpha\tau^2\right)^{-{1 \over 2\alpha}}=\exp\left\{ -{1 \over 2\alpha} \ln \left(1+\alpha \tau^2\right) \right\}.$$
Then
$${\partial \over \partial \alpha} h(\tau, \alpha) =\left({1 \over 2 \alpha^2} \ln \left(1+ \alpha \tau^2\right)  -{\tau^2 \over 2 \alpha \left(1+ \alpha \tau^2\right)}\right)\exp\left\{ -{1 \over 2 \alpha} \ln \left( 1 + \alpha\tau^2\right)\right\}$$
and
$${1 \over 2 \alpha^2} \ln \left(1+ \alpha \tau^2\right)  -{\tau^2 \over 2 \alpha \left(1+ \alpha \tau^2\right)}=  {(1+\alpha \tau^2)  \ln \left(1+ \alpha \tau^2\right) - \alpha \tau^2 \over 2 \alpha^2 (1+\alpha \tau^2)}.$$
Finally, we observe that 
$$g(\sigma) = (1+\sigma) \ln (1+ \sigma) - \sigma \ > \ 0 \quad \text{for} \quad \sigma > 0,$$
since $g(0)=0$ and 
$$g'(\sigma)=\ln (1+ \sigma) \ > \ 0 \quad \text{for} \quad \sigma > 0.$$
Summarizing, the function $\alpha \longmapsto h(\tau, \alpha)$ is increasing for all $\tau \ne 0$ and constant for $\tau=0$. The proof of Part (1) follows.

For $0 < \alpha < 1$, we have 
$$\sqrt{\alpha} F(\alpha)={\sqrt{\alpha} \over 2\pi} \int_{-\infty}^{+\infty} \left(1 + \alpha \tau^2\right)^{-{1 \over 2\alpha}} \ d \tau = {1 \over 2 \pi} \int_{-\infty}^{+\infty} \left(1+ \sigma^2\right)^{-{1 \over 2\alpha}} \ d \sigma.$$
Hence as $\alpha$ changes from $0$ to $1$, the value of $\sqrt{\alpha}F(\alpha)$ increases from $0$ to $+\infty$. We find $\alpha_0$ from the equation $\sqrt{\alpha_0}F(\alpha_0)=1$,
which we solve numerically. This completes the proof of Part (2).
{\hfill \hfill \hfill} \qed
\enddemo

Recall that $b_1, \ldots, b_n$ are the columns of matrix $B$. By $\|\cdot \|$ we denote the standard Euclidean norm in ${\Bbb R}^m$. 
\proclaim{(4.4) Corollary} Suppose that the columns $b_1, \ldots, b_n$ of matrix $B$ satisfy
$$\sum_{j=1}^n b_j \otimes b_j =I_m $$ and that 
$$\|b_j\| \ \leq \ \sqrt{\alpha_0} \quad \text{for} \quad j=1, \ldots, n,$$
where $\alpha_0$ is the constant of Lemma 4.3. Then 
$$\left\| p_Y\right\|_{\infty} \ \leq \ \alpha_0^{-{m/2}},$$
where $Y$ is the random variable defined by (3.1.1).
\endproclaim

\demo{Proof} We use (4.1.1).
Without loss of generality, we assume that $b_j \ne 0$ for $j=1, \ldots, n$. Let 
$$\lambda_j = \|b_j\|^2 \quad \text{and} \quad u_j = {1 \over \sqrt{\lambda_j}} b_j \quad \text{for} \quad j=1, \ldots, n.$$
Hence $u_j$ are unit vectors,
$$\sum_{j=1}^n \lambda_j \left(u_j \otimes u_j \right) = I_m \quad \text{and} \quad \sum_{j=1}^n \lambda_j=m, \tag4.4.1$$
where the second identity is obtained comparing the traces of matrices on both sides of the first identity.
Besides.
$$0 < \lambda_j \leq \alpha_0 \quad \text{for} \quad j=1, \dots, n.$$
By (4.1.1), we have 
$$\split \left\| p_Y \right\|_{\infty} \ \leq \ &{1 \over (2 \pi)^m} \int_{{\Bbb R}^m} \prod_{j=1}^n \left( 1+ \lambda_j \langle u_j, t \rangle^2 \right)^{-1/2} \ dt \\
= & {1 \over (2 \pi)^m} \int_{{\Bbb R}^m} \prod_{j=1}^n \left(\left( 1+ \lambda_j \langle u_j, t \rangle^2 \right)^{-{1 \over 2 \lambda_j}}\right)^{\lambda_j} \ dt  \\
\ \leq \ & {1 \over (2 \pi)^m} \prod_{j=1}^n \left( \int_{-\infty}^{+\infty} \left(1 + \lambda_j \tau^2 \right)^{-{1 \over 2\lambda_j}} \ d \tau\right)^{\lambda_j}  \\ 
\ \leq \ & {1 \over (2 \pi)^m} \ \prod_{j=1}^n \left( \int_{-\infty}^{+\infty} \left(1 + \alpha_0 \tau^2 \right)^{-{1 \over 2\alpha_0}} \ d \tau\right)^{\lambda_j} \\
\ = \ &{1 \over (2 \pi)^m} \prod_{j=1}^n \left( {2 \pi \over \sqrt{\alpha_0}}\right)^{\lambda_j} =\alpha_0^{-m/2}.  \endsplit $$
We use Lemma 4.2 in the inequality of the third line, Part (1) of Lemma 4.3 in the inequality of the fourth line and Part (2) of Lemma 4.3 and (4.4.1) in the last line.
{\hfill \hfill \hfill} \qed
\enddemo

To complete the proof, we need the following standard result.
\proclaim{(4.5) Lemma} Let $\varphi, \psi: {\Bbb R} \longrightarrow {\Bbb R}_+$ be densities and let 
$\rho: {\Bbb R} \longrightarrow {\Bbb R}_+$ be their convolution 
$$\rho(\xi) = \int_{-\infty}^{+\infty} \varphi(\xi-\tau) \psi(\tau) \ d \tau.$$
Then $\rho$ is a density and 
$$ \| \rho\|_{\infty} \ \leq \ \min \left\{ \| \varphi\|_{\infty}, \ \|\psi \|_{\infty} \right\}.$$
\endproclaim
{\hfill \hfill \hfill} \qed 

\subhead (4.6) Proof of Part (1) of Theorem 2.2 \endsubhead 
As before, we assume that the matrix $B$ satisfies (3.4.2), or equivalently, the columns $b_1, \ldots, b_n$ of $B$ satisfy 
$$\sum_{j=1}^n b_j \otimes b_j = I_m, \tag4.6.1$$
where $I_m$ is the $m \times m$ identity matrix. Another equivalent way to write (4.6.1) is 
$$\sum_{j=1}^n \langle b_j, x \rangle^2 = \|x\|^2 \quad \text{for all} \quad x \in {\Bbb R}^m.$$
 In view of (3.4.1), our goal is to extend the conclusion of Corollary 4.4, without assuming that 
$\|b_j\| \leq \sqrt{\alpha_0}$ for $j=1, \ldots, n$.

Let $X_1, \ldots, X_n$ be independent standard exponential random variables and let $Y$ be defined by (3.1.1).
We proceed by induction on $m$. Suppose that $m=1$, so
$$Y= \sum_{j=1}^n \mu_j X_j \quad \text{where} \quad \sum_{j=1}^n \mu_j^2 =1.$$
If we have 
$$|\mu_j| \ \leq \ \sqrt{\alpha_0} \quad \text{for} \quad j=1, \ldots, n,$$
the result follows from Corollary 4.4. If for at least one $\mu_j$ we have $|\mu_j |> \sqrt{\alpha_0}$ and for some $i \ne j$ we have $\mu_i \ne 0$ then applying Lemma 4.5, where
$\varphi=p_{\mu_j X_j}$, $\psi$ is the density of $\sum_{i:\ i \ne j} \mu_i X_i$, and $\rho=p_Y$,
 we obtain 
$$\left\| p_Y\right\|_{\infty} \ \leq \ \left\| p_{\mu_j X_j} \right\|_{\infty} = {1 \over |\mu_j|} \ < \ {1 \over \sqrt{\alpha_0}}.$$
If $\mu_i =0$ for all $i \ne j$ then $|\mu_j|=1$ and 
$$\left\| p_Y\right\|_{\infty} = \left\| p_{X_j}\right\|_{\infty} = 1 < {1 \over \sqrt{\alpha_0}}.$$

Suppose now that $m >1$. Let 
$$\lambda_j = \|b_j\|^2 \quad \text{for} \quad j=1, \ldots, n.$$
If 
$$\lambda_j \ \leq \ \alpha_0 \quad \text{for} \quad j=1, \ldots, n, $$
the result follows by Corollary 4.4. Otherwise, we have $\lambda_j > \alpha_0$ for some $j$. Without loss of generality, we assume that 
$\lambda_n > \alpha_0$. From (4.6.1) it follows that $\lambda_n \leq 1$.

Suppose first that $\lambda_n=1$. Then from (4.6.1) we must have
$$\langle b_j, b_n \rangle =0 \quad \text{for} \quad j=1, \ldots, n-1. \tag4.6.2$$
We consider a decomposition ${\Bbb R}^m ={\Bbb R}^{m-1} \oplus {\Bbb R}$, $y=(y', \eta)$,
where we identify 
$$\spa\left(b_1, \ldots, b_{n-1}\right)= {\Bbb R}^{m-1} \quad \text{and} \quad  \spa\left(b_n\right) = {\Bbb R}.$$
From (4.6.1) and (4.6.2), we have 
$$\sum_{j=1}^{n-1} b_j \otimes b_j = I_{m-1}.$$
Let $Y'=\sum_{j=1}^{n-1} X_j b_j $
be a random vector in ${\Bbb R}^{m-1}$ and let $Y'' = X_n b_n = \pm X_n$ be a random variable with values in ${\Bbb R}$, so that $Y=\left(Y', Y''\right)$. 
Then 
$$p_Y(y', \eta)=p_{Y'} (y') p_{Y''}(\eta)$$ 
and applying the induction hypothesis, we obtain 
$$\left\| p_Y \right\|_{\infty} = \left\| p_{Y'}\right\|_{\infty} \left\| p_{Y''} \right\|_{\infty} = \left\|p_{Y'}\right\|_{\infty} \ \leq \ \alpha_0^{-(m-1)/2} \ < \ \alpha_0^{-m/2}.$$
It remains to consider the case where 
$$ \alpha_0 \ < \ \lambda_n \ < \ 1. \tag4.6.3$$
We consider a decomposition ${\Bbb R}^m ={\Bbb R}^{m-1} \oplus {\Bbb R}$, where ${\Bbb R}$ is identified with $\spa\left(b_n\right)$ and 
${\Bbb R}^{m-1}$ is identified with the orthogonal complement $b_n^{\bot}$. For $j=1, \ldots, n-1$, let $b_j'$ be the orthogonal projection of $b_j$ onto ${\Bbb R}^{m-1}$
and let $b_j''$ be the orthogonal projection of $b_j$ onto ${\Bbb R}$. From (4.6.1) it follows that 
$$\sum_{j=1}^{n-1} b_j' \otimes b_j'= I_{m-1}.$$
We introduce a random vector $Y'$ with values in ${\Bbb R}^{m-1}$ by
$$Y'=\sum_{j=1}^{n-1} X_j b_j'$$
and a random variable $Y''$ with values in ${\Bbb R}$ by 
$$Y'' =X_n b_n + \sum_{j=1}^{n-1} X_j b_j'',$$
so that $Y=(Y', Y'')$.
By the induction hypothesis, we have 
$$\left\| p_{Y'}\right\|_{\infty} \ \leq \ \alpha_0^{-(m-1)/2}. \tag4.6.4$$
Using conditional density, for $y \in {\Bbb R}^m$, $y=(y', \eta)$ where $y' \in {\Bbb R}^{m-1}$ and $\eta \in {\Bbb R}$, we write 
$$p_Y(y', \eta)= \cases 0 &\text{if\ } p_{Y'}(y')=0  \\ p_{Y'}(y') p_{Y''|Y'}(\eta | y') &\text{if\ } p_{Y'}(y') \ne 0. \endcases \tag4.6.5$$
Let
$$Z=\sum_{j=1}^{n-1} X_j b_j''.$$
Since (4.6.1) and (4.6.3) hold, the matrix $\sum_{i=1}^{n-1} b_i \otimes b_i$ is invertible.  Therefore, the random vector
$$(Y', Z) = \sum_{j=1}^{n-1} X_j b_j$$
has density and hence the conditional density $p_{Z|Y'}$ exists whenever $p_{Y'} \ne 0$.

Now, we have $Y''=Z + X_n b_n$ and hence
for the conditional densities we have 
$$p_{Y''|Y'}= p_{Z| Y'} \ast p_{X_n b_n}.$$
Applying Lemma 4.5 with $\rho=p_{Y''|Y'}$, $\varphi=p_{X_n b_n}$ and $\psi=p_{Z|Y'}$, we obtain
$$\left\| p_{Y'' | Y'} \right\|_{\infty} \ \leq \ \left\| p_{X_n b_n}\right\|_{\infty} ={1 \over \|b_n\|} \ \leq \ {1 \over \sqrt{\alpha_0}}. \tag4.6.6$$
Combining (4.6.4) -- (4.6.6), we conclude that 
$$\left\| p_Y\right\|_{\infty} \ \leq \ \alpha_0^{-m/2}.$$
The proof now follows from (3.4.1).
{\hfill \hfill \hfill} \qed

\head 5. Proof of Part (2) \endhead 

The bound of Part (2) of Theorem 2.2 will follow from some general estimate for isotropic log-concave densities.
\proclaim{(5.1) Lemma} Let $\varphi: {\Bbb R}^m  \longrightarrow {\Bbb R}_+$ be a centered log-concave density and let $H \subset {\Bbb R}^n$ be a closed halfspace containing $0$. Then 
$$\int_H \varphi(x) \ dx \ \geq \ {1 \over e}.$$
\endproclaim

\demo{Proof} When $\varphi$ is the uniform density on a convex body, the result was proved by Gr\"unbaum \cite{Gr60} (with a slightly better constant depending on $n$ and decreasing to $1/e$ as $n$ grows). For an adaptation to general log-concave measures, see \cite{LV07}, Proposition 1.5.16 in \cite{A+15} or Lemma 2.2.6 in \cite{B+14}.
{\hfill \hfill \hfill} \qed
\enddemo

\proclaim{(5.2) Theorem} Let $\varphi: {\Bbb R}^m  \longrightarrow {\Bbb R}_+$ be an isotropic log-concave density. Then 
$$\varphi(0) \ \geq \   {2 \Gamma\left({m +2 \over 2}\right) \over \pi^{m/2} e^{(m+2)/2} (m+2)^{m/2}}.$$
\endproclaim
\demo{Proof} Without loss of generality, we assume that $\varphi$ is not constant on open subsets of ${\Bbb R}^n$: a general $\varphi$ can be approximated by a sequence of strictly log-concave isotropic densities $\varphi_n$ obtained from $\varphi(x) \exp\left\{ - \|x\|^2/n\right\}$ by a scaling, shift and linear transformation. Then the set
$$K=\Bigl\{ x \in {\Bbb R}^m: \ \varphi(x) \geq \varphi(0) \Bigr\}$$
is convex and such that $0 \in \partial K$. There is a hyperplane supporting $K$ at $0$. Let $H$ be the open halfspace bounded by that hyperplane and disjoint from the interior of $K$. We have 
$$\varphi(x) \ \leq \ \varphi(0) \quad \text{for all} \quad x \in H \tag5.2.1$$
and by Lemma 5.1, 
$$\beta:=\int_H \varphi(x)\ dx \ \geq \ {1 \over e}. \tag5.2.2$$ 
It is clear now that $\varphi(0) >0$.
From this point on, our proof mimics that of Proposition 10.2.5 of \cite{A+15} that provides a lower bound for the $\ell^{\infty}$-norm of an isotropic, but not necessarily log-concave density. 
Let
$$\kappa_m={\pi^{m/2} \over \Gamma\left({m+2 \over 2}\right)}$$
denote the volume of the unit ball in ${\Bbb R}^m$ and let 
$$D_{\tau}=\Bigl\{x \in \ {\Bbb R}^m: \ \|x\| \leq \tau \Bigr\}$$
denote the ball of radius $\tau >0$ in ${\Bbb R}^m$.
Since $\varphi$ is isotropic, we have
$$\int_{{\Bbb R}^m} \|x\|^2 \varphi(x) \ dx =\sum_{i=1}^m \int_{{\Bbb R}^m} \xi_i^2 \varphi(x) \ dx =m. \tag5.2.3$$
Let $\rho >0$ be a number, to be specified later.

Taking into account (5.2.3) and then (5.2.1), we obtain
 $$\split m = &\int_{{\Bbb R}^m} \|x\|^2 \varphi(x) \ d x \ \geq \ \int_{H} \|x\|^2 \varphi(x) \ dx  =\int_{H} \left( \int_0^{\|x\|^2} 1 \ d \tau \right) \varphi(x) \ dx \\
 =&\int_0^{+\infty} \left(\int_{H \setminus D_{\sqrt{\tau}}} \varphi(x) \ dx \right) \ d\tau = \int_0^{+\infty} \left(\beta -\int_{H
 \cap  D_{\sqrt{\tau}}} \varphi(x) \ dx \right) \ d \tau \\
 \ \geq & \int_0^{\rho} \left(\beta -\int_{H \cap  D_{\sqrt{\tau}}} \varphi(x) \ dx \right) \ d \tau 
  \ \geq \ \int_0^{\rho} \left(\beta - {1 \over 2} \kappa_m \varphi(0) \tau^{m/2} \right) \ d\tau \\
  =\ & \rho \beta - {\kappa_m f(0) \over m+2} \rho^{m+2 \over 2}.\endsplit $$
  Optimizing on $\rho$, we choose
  $$\rho=\left({2\beta \over \kappa_m \varphi(0)}\right)^{2/m} $$
  and obtain 
  $$\split m \ \geq \ &\left({2 \beta \over \kappa_m \varphi(0)}\right)^{2/m} \beta - {\kappa_m \varphi(0) \over m+2} \left({2 \beta \over \kappa_m \varphi(0)}\right)^{m+2 \over m}\\
  =& {2^{2/m} \beta^{(m+2)/m} \over \left( \kappa_m \varphi(0)\right)^{2/m}}  {m \over m+2},\endsplit$$
 from which
 $$\varphi(0) \ \geq \ {2 \beta^{m+2 \over 2}  \over \kappa_m (m+2)^{m/2}}= {2 \beta^{m+2 \over 2} \Gamma\left({m+2 \over2}\right) \over \pi^{m/2}  (m+2)^{m/2}}, $$
 and the proof follows by (5.2.2).
{\hfill \hfill \hfill} \qed
\enddemo

\subhead (5.3) Proof of Part (2) \endsubhead As before, without loss of generality, we assume that matrix $B$ satisfies (3.4.2). Then by (3.1.3), the density $p_{Y-b}$ of $Y-b$ 
is isotropic. It is also, as we discussed in Section 3.4, log-concave. Hence by Theorem 5.2, we have 
$$p_Y(b) = p_{Y-b}(0) \ \geq \ {2 \Gamma \left({m +2 \over 2}\right) \over \pi^{m/2} e^{(m+2)/2} (m+2)^{m/2}}.$$
The proof now follows by (3.4.1).
{\hfill \hfill \hfill} \qed

\head 6. Proof of Part (3) \endhead

\subhead (6.1) More preliminaries: thin shell estimates \endsubhead Let $\varphi: {\Bbb R}^m \longrightarrow {\Bbb R}$ be an isotropic log-concave density. Then
$$m=\int_{{\Bbb R}^m} \|x\|^2 \varphi(x) \ dx$$ is the expectation of $\|x\|^2$ with respect to the density $\varphi$. In a recent paper, 
Klartag and and Lehec \cite{KL22} established the following inequality, known as a ``thin-shell estimate", see Chapter II of \cite{A+21}, for the variance of $\|x\|^2$:
$$\int_{{\Bbb R}^m} \left(\|x\|^2-m \right)^2 \varphi(x) \ dx  \ \leq \ \gamma m \ln^8 m \tag6.1.1$$
for some absolute constant $\gamma > 0$ and all $m > 1$.

Let $\mu$ be the probability measure on ${\Bbb R}^m$ with density $\varphi$. The bound (6.1.1) implies that the value of $\|x\|$ is strongly concentrated about $\sqrt{m}$. In particular, 
by the Chebyshev inequality, for any $\tau >0$, we get 
$$\aligned &\mu\Bigl\{ x \in {\Bbb R}^m:\ \|x\| \geq \sqrt{m} + \tau \Bigr\} = \mu\left\{ x \in {\Bbb R}^m:\ \|x\|^2 \ \geq \ m + 2 \tau \sqrt{m} + \tau^2 \right\} \\ &\qquad \leq \ 
{\gamma m \ln^8 m \over (2 \tau \sqrt{m} + \tau^2)^2}. \endaligned \tag6.1.2$$
To obtain the main asymptotic 
$$\left({1 \over \sqrt{2 \pi e}}\right)^m$$ in Part (3) of Theorem 2.2 and Theorem 6.2 below, we could have used the earlier thin shell bounds by 
Chen \cite{Ch21}
$$\mu\Bigl\{x \in {\Bbb R}^m: \ \|x\| \geq \sqrt{m} + \tau \Bigr\} \ \leq \ \tau^{-1} \exp\left\{\gamma_1 \sqrt{\ln m \ln \ln m}\right\}$$
for $\tau >0$, absolute constant $\gamma_1 >0$ and all $m \geq 3$ and Klartag \cite{Kl07}
$$\mu\Bigl\{x \in {\Bbb R}^m: \ \|x\| \geq \sqrt{m} (1+ \epsilon ) \Bigr\} \ \leq \ \gamma_2 m^{-\gamma_3 \epsilon^2}$$
for $0 \leq \epsilon \leq 1$ and absolute constants $\gamma_2, \gamma_3 >0$, see also \cite{LV17} for a survey.

\proclaim{(6.2) Theorem}  For any $0 < \epsilon < 1/2$, there is $\gamma(\epsilon)>0$ such that if $\varphi: {\Bbb R}^m \longrightarrow {\Bbb R}_+$ is an isotropic log-concave density, then 
$$\varphi(0) \ \geq \ \exp\left\{{- \gamma(\epsilon)} \sqrt{m}  \ln^4 (m+1) \right\} \left( {1 \over \sqrt{2 \pi e}}\right)^m.$$
\endproclaim
\demo{Proof} In view of Theorem 5.2, without loss of generality, we assume that $m> 1$.  As in the proof of Theorem 5.2, we find a halfspace $H \subset {\Bbb R}^m$
such that $0 \in \partial H$ so that 
$$\int_{H} \varphi(x) \ dx \ \geq \ {1 \over e} \quad \text{and} \quad \varphi(0) \ \geq \ \varphi(x) \quad \text{for all} \quad  x\in H. \tag6.2.1$$
Let $\mu$ be the probability measure on ${\Bbb R}^m$ with density $\varphi$. Choosing 
$$\tau = \sqrt{2 \gamma} \ln^4 m$$
in (6.1.2), 
we conclude that 
$$\mu \Bigl\{ x \in {\Bbb R}^m: \ \|x\| \geq \ \sqrt{m} + \sqrt{2 \gamma} \ln^4 m \Bigr\} \ \leq \ {1 \over 8} \ < \ {1 \over 2 e}$$
and hence by (6.2.1) we have 
$$\mu\Bigl\{ x \in H:\  \|x \| \leq \ \sqrt{m} +  \sqrt{2 \gamma} \ln^4 m \Bigr\} \ \geq \ {1 \over 2e}.$$
Therefore,
$$\varphi(0) \kappa_m \left(\sqrt{m} + \sqrt{2 \gamma} \ln^4 m \right)^m \ \geq \ {1 \over e},$$
where $\kappa_m$ is the volume of the unit ball in ${\Bbb R}^m$, and 
$$\phi(0) \ \geq \ {1 \over  e \kappa_m} \left(\sqrt{m} + \sqrt{2 \gamma} \ln^4 m \right)^{-m} ={\Gamma\left({m+2 \over 2} \right) \over e \pi^{m/2}
 \left( \sqrt{m} + \sqrt{2 \gamma} \ln^4 m  \right)^{m}}. $$
The proof now follows by (2.2.1).
{\hfill \hfill \hfill} \qed
\enddemo 

\subhead (6.3) Proof of Part (3) \endsubhead The proof follows as in Section 5.3, except that we use Theorem 6.2 instead of Theorem 5.2.
{\hfill \hfill \hfill} \qed

\enddocument

\end